\documentclass{amsart}

\usepackage{amssymb, amscd, latexsym, graphicx, psfrag}

\newcommand{\Ga}{\Gamma}

\newcommand{\bC}{\mathbb{C}}

\newcommand{\bQ}{\mathbb{Q}}

\newcommand{\bZ}{\mathbb{Z}}

\newcommand{\cB}{\mathcal{B}}
\newcommand{\cC}{\mathcal{C}}
\newcommand{\cD}{\mathcal{D}}
\newcommand{\cE}{\mathcal{E}}
\newcommand{\cF}{\mathcal{F}}
\newcommand{\cI}{\mathcal{I}}
\newcommand{\cL}{\mathcal{L}}
\newcommand{\cM}{\mathcal{M}}

\newcommand{\cR}{\mathcal{R}}
\newcommand{\cQ}{\mathcal{Q}}

\newcommand{\cX}{\mathcal{X}}

\newcommand{\Hom}{\mathrm{Hom}}

\newcommand{\age}{\mathrm{age}}

\newcommand{{\inv} }{\mathrm{inv}}
\newcommand{\ev}{\mathrm{ev}}
\newcommand{\Aut}{\mathrm{Aut}}

\newcommand{\rank}{\mathrm{rank}}
\newcommand{\val}{ {\mathrm{val}} }
\newcommand{\tw}{ {\mathrm{tw}} }
\newcommand{\CR}{ {\mathrm{CR}} }

\newcommand{\zero}{\mathbf{0}}
\newcommand{\one}{\mathbf{1}}

\newcommand{\be}{\mathbf{e}}
\newcommand{\bu}{\mathbf{u}}

\newcommand{\bGa}{\mathbf{\Ga}}

\newcommand{\su}{\mathsf{u}}

\newcommand{\sw}{\mathsf{w}}

\newcommand{\tcL}{\widetilde{\cL} }
\newcommand{\tT}{ \widetilde{\mathbb{T}} }

\newcommand{\hbu}{\hat{\bu}}

\newcommand{\Mbar}{\overline{\cM}}

\newcommand{\vGa}{\vec{\Ga}}

\newcommand{\BG}{\cB G}
\newcommand{\IBG}{\cI \cB G}
\newcommand{\IX}{\cI \cX}

\newtheorem{dummy}{dummy}

\newtheorem{theorem}[dummy]{Theorem}

\theoremstyle{definition}

\newtheorem*{remark*}{Remark}

\begin{document}

\title[Equivariant GW Theory of Affine Smooth Toric DM Stacks]{Equivariant
Gromov-Witten Theory of Affine Smooth Toric Deligne-Mumford Stacks}

\author{Bohan Fang}
\address{Bohan Fang, Beijing International Center for Mathematical Research, Peking University, 5 Yiheyuan Road, Beijing 100871, China}
\email{bohanfang@gmail.com}

\author{Chiu-Chu Melissa Liu}
\address{Chiu-Chu Melissa Liu, Department of Mathematics, Columbia University, 2990 Broadway, New York, NY 10027}
\email{ccliu@math.columbia.edu}

\author{Zhengyu Zong}
\address{Zhengyu Zong, Department of Mathematics, Columbia University, 2990 Broadway, New York, NY 10027}
\email{zz2197@math.columbia.edu}

\maketitle

\begin{abstract}
For any finite abelian group $G$, the equivariant Gromov-Witten invariants of $[\bC^r/G]$ can be viewed as a certain kind of abelian Hurwitz-Hodge integrals. In this note, we use Tseng's orbifold quantum Riemann-Roch theorem to express this kind of abelian Hurwitz-Hodge integrals as a sum over Feynman graphs, where the weight of each graph is expressed in terms of descendant integrals over moduli spaces of stable
curves and representations of $G$. This expression will play a crucial role in the proof of the remodeling conjecture for affine toric Calabi-Yau 3-orbifolds in \cite{FLZ}.

\end{abstract}

\section{Introduction}

Gromov-Witten invariants are virtual counts of parametrized algebraic curves in algebraic manifolds, or more generally,
parametrized holomorphic curves in K\"{a}hler manifolds. Gromov-Witten theory can be viewed
as a mathematical theory of A-model topological string theory. The topological string theory on orbifolds was constructed 
decades ago by physicists \cite{DHVW85, DHVW86}, and many works followed in both mathematics and physics (e.g. \cite{HV, DFMS, BR, Kn, Roan, CV}).  Orbifold quantum cohomology was discussed in \cite{Za} along with many examples, in both abelian and non-abelian quotients. Later, the mathematical definition of orbifold Gromov-Witten theory and quantum cohomology was laid in \cite{CR02, AGV02, AGV08}.

Let $\tT=(\bC^*)^r$ act on $\bC^r$ by $(t_1,\ldots, t_r)\cdot (z_1,\ldots, z_r)=(t_1 z_1,\ldots, t_r z_r)$,
and let $\phi: G\to \tT$ be a group homomorphism, where $G$ is a finite abelian group.
The quotient stack $[\bC^r/G]$ is an affine smooth toric DM stack (in the sense of \cite{BCS05}); 
it is an affine toric orbifold iff $\phi$ is injective.
In this note, we study all genus descendant $\tT$-equivariant
orbifold Gromov-Witten invariants of $[\bC^r/G]$, which are essentially
abelian Hurwitz-Hodge integrals. By virtual localization,
orbifold Gromov-Witten theory of a general smooth toric DM stack can
be reduced to theory of abelian Hurwitz-Hodge integrals.

The quotient stack $[\bC^r/G]$ can be viewed as the total space of a
$\tT$-equivariant rank $r$ vector bundle over the classifying space
$\BG$ of $G$, where $\tT$ acts trivially on the base $\BG$.
By Tseng's orbifold quantum Riemann-Roch theorem \cite{Ts10},
all genus descendant $\tT$-equivariant orbifold Gromov-Witten invariants of $[\bC^r/G]$
can be reconstructed from descendant integrals over moduli spaces of
twisted stable maps to $\BG$. By Jarvis-Kimura \cite{JK}, descendant integrals over moduli spaces
of twisted stable maps to $\BG$ can be reduced to descendant integrals over moduli spaces of
stable curves $\Mbar_{g,n}$.  The goal of this note is to make the above process more explicit
in terms of Feynman graphs, following \cite{DSS, DOSS}. Our main result Theorem \ref{main}
is an expression of the total descendant $\tT$-equivariant Gromov-Witten potential of 
$[\bC^r/G]$ as a sum over graphs, where the weights assigned
to vertices of these graphs are essentially descendant integrals over $\Mbar_{g,n}$.

Theorem \ref{ordered} is a generalized version of
Theorem \ref{main}.  (See Section \ref{sec:graph-sum} for the precise statements
of Theorem \ref{main} and Theorem \ref{ordered}.) In \cite{FLZ}, we use Theorem \ref{ordered} to prove
the remodeling conjecture for affine toric
Calabi-Yau 3-orbifolds $[\bC^3/G]$, where $G$ is a finite
subgroup of the maximal torus of $SL(3,\bC)$. The remodeling conjecture proposed
by Bouchard-Klemm-Mari\~{n}o-Vafa \cite{BKMP09, BKMP10}
relates all genus open and closed Gromov-Witten invariants
of a toric Calabi-Yau 3-orbifold $\cX$ (in particular, a smooth
Calabi-Yau 3-folds) to Eynard-Orantin invariants \cite{EO07} of the mirror
curve of $\cX$. The remodeling conjecture was proved for $\bC^3$ independently
by L. Chen \cite{Ch09} and by J. Zhou \cite{Zh09}. In \cite{EO12}, 
Eynard-Orantin provide a proof of the remodeling conjecture 
for general smooth toric Calabi-Yau 3-folds.

\subsection*{Acknowledgments} The second author wishes to thank B. Eynard and N. Orantin for helpful conversations.
The research of the first author is partially supported by NSF DMS-1206667.
The research of the second and third authors is partially supported by NSF DMS-1159416.

\section{Chen-Ruan Orbifold Cohomology}

In this section, we describe the Chen-Ruan orbifold cohomology of $\BG$ and the $\tT$-equivariant Chen-Ruan orbifold cohomology of $[\bC^r/G]$ \cite{CR04}. The quantum cohomology of an orbifold
is a deformation of its Chen-Ruan orbifold cohomology. The consideration of inertia stack in orbifolds 
appears very naturally in the study of orbifolds, like orbifold Riemann-Roch theorem \cite{Ka}.

\subsection{Chen-Ruan orbifold cohomology of $\BG$}
Let $G$ be a finite abelian group and $\BG$ the classifying stack of $G$. The inertia stack of $\BG$ is
$$
\IBG = \bigcup_{h\in G}(\BG)_h,
$$
where
$$
(\BG)_h =[\{h\}/G]\cong \BG.
$$
As a graded vector space over $\bC$, the Chen-Ruan orbifold
cohomology of $\BG$ is
$$
H^*_\CR(\BG;\bC) = H^*(\IBG;\bC) = \bigoplus_{h\in G} H^0((\BG)_h;\bC),
$$
where $H^0(\BG_h;\bC)=\bC\one_h$.
The orbifold Poincar\'{e} pairing of $H^*_\CR(\BG;\bC)$ is given by
$$
\langle\one_h, \one_{h'}\rangle = \frac{\delta_{h^{-1},h'}}{|G|}.
$$
The orbifold cup product of $H^*_\CR(\BG;\bC)$ is given by
$$
\one_h \star \one_{h'} = \one_{h h'}.
$$

Following \cite{JK}, we define a canonical basis for
the semisimple algebra $H^*_\CR(\BG;\bC)$.
Given a character $\gamma\in G^*=\Hom(G,\bC^*)$, define
$$
\phi_\gamma :=\frac{1}{|G|}\sum_{h\in G}\chi_\gamma(h^{-1})\one_h.
$$
Then
$$
H^*_\CR(\BG;\bC) =\bigoplus_{h\in G}\bC \one_h  =\bigoplus_{\gamma \in G^*}\bC \phi_\gamma.
$$
Recall that we have the orthogonality of characters:
\begin{enumerate}
\item For any $\gamma,\gamma'\in G^*$,
$\displaystyle{
\frac{1}{|G|} \sum_{h\in G} \chi_{\gamma}(h^{-1})\chi_{\gamma'}(h)=\delta_{\gamma,\gamma'}
},$
\item For any $h, h'\in G$,
$\displaystyle{
\frac{1}{|G|}\sum_{\gamma\in G^*}\chi_\gamma(h^{-1})\chi_\gamma(h') =\delta_{h,h'}.
}$
\end{enumerate}
Therefore,
$$
\langle \phi_{\gamma}, \phi_{\gamma'}\rangle = \frac{\delta_{\gamma,\gamma'}}{|G|^2},
$$
and
$$
\phi_{\gamma}\star\phi_{\gamma'}= \delta_{\gamma, \gamma'}\phi_{\gamma}.
$$
In particular, $\{\phi_\gamma: \gamma\in G^*\}$ is a canonical basis
of the semisimple $\bC$-algebra $H^*_\CR(\cB G;\bC)$.

\subsection{Equivariant Chen-Ruan orbifold cohomology of $\cX=[\bC^r/G]$}
Given any $h\in G$ and $i\in \{1,\ldots, r\}$, define
$c_i(h)$ by
$$
\chi_i(h) = e^{2\pi\sqrt{-1} c_i(h)},\quad 0\leq c_i(h)<1.
$$
where $\chi_i: G\to \bC^*$ is the character associated to the action of $G$ on the $i$-th component of $\bC^r$. And define the age of $h$ to be
$$
\age(h)=\sum_{i=1}^n c_i(h).
$$
Let $(\bC^r)^h$ denote the $h$-invariant subspace of $\bC^r$. Then
$$
\dim_\bC (\bC^r)^h  = \sum_{i=1}^r \delta_{c_i(h),0}.
$$
The inertial stack of $\cX$ is
$$
\IX  = \bigcup_{h\in G}\cX_h,
$$
where
$$
\cX_h =[(\bC^r)^h/G].
$$
In particular,
$$
\cX_1 =[\bC^r/G]=\cX.
$$
As a graded vector space over $\bC$,
$$
H^*_\CR(\cX;\bC) = \bigoplus_{h\in G} H^*(\cX_h; \bC)[2\age(h)] =\bigoplus_{h\in G} \bC \one_h,
$$
where $\deg(\one_h)=2\age(h)\in \bQ$. The orbifold Poincar\'{e} pairing of the (non-equivariant) Chen-Ruan orbifold cohomology $H_\CR^*(\cX;\bC)$ is given by
$$
\langle \one_{h}, \one_{h'}\rangle_{\cX} =\frac{1}{|G|}\delta_{h^{-1},h'}
\cdot \delta_{0,\dim_\bC \cX_h}.
$$

Let $\cR=H^*(B\tT;\bC)=\bC[\sw_1,\ldots,\sw_r]$, where
$\sw_1,\ldots,\sw_r$ are the first Chern classes of the universal line bundles over $B\tT$.
The $\tT$-equivariant Chen-Ruan orbifold cohomology $H^*_{\CR,\tT}([\bC^r/G];\bC)$
is an $\cR$-module.  Given $h\in G$, define
$$
\be_h := \prod_{i=1}^r \sw_i^{\delta_{c_i(h),0}} \in \cR.
$$
In particular,
$$
\be_1 =\prod_{i=1}^r \sw_i.
$$
Then the $\tT$-equivariant Euler class of $\zero_h:=[0/G]$ in $\cX_h = [(\bC^r)^h/G]$ is
$$
e_{\tT}(T_{\zero_h}\cX_h) = \be_h \one_h \in H^*_{\tT}(\cX_h;\bC) = \cR \one_h.
$$

Recall that the $G$-action on the $i$-th axis of $\bC^r$ determines a character
$\chi_i \in G^*=\Hom(G,\bC^*)$. The image of $\chi_i:G\to \bC^*$ is
isomorphic to $\bZ/l_i\bZ$ for some positive integer $l_i$. Define
$$
\cR' =\bC[\sw_1^{1/l_1},\ldots, \sw_r^{1/l_r}],
$$
which is a finite extension of $\cR$. Let $\cQ$ and $\cQ'$ be the fractional fields of 
$\cR$ and $\cR'$, respectively. The $\tT$-equivariant Poincar\'{e} pairing of
$H^*_{\tT,\CR}(\cX;\bC)\otimes_{\cR}\cQ$ (which is isomorphic to  $H^*_{\CR}(\cX;\cQ)$
as a vector space over $\cQ'$) is given by
$$
\langle \one_{h}, \one_{h'}\rangle_{\cX} =\frac{1}{|G|}\cdot\frac{\delta_{h^{-1},h'} }{\be_h}\in \cQ.
$$
The $\tT$-equivariant orbifold cup product  of $H^*_{\tT,\CR}(\cX;\bC)\otimes_{\cR}\cQ$ is given by
$$
\one_h \star_{\cX} \one_{h'} = \Bigl(\prod_{i=1}^r \sw_i^{c_i(h)+c_i(h')-c_i(hh')}\Bigr) \one_{hh'}.
$$

Define
\begin{equation}\label{eqn:bar-one}
\bar{\one}_h:= \frac{\one_h}{\prod_{i=1}^r \sw_i^{c_i(h)}} \in H^*_{\tT,\CR}(\cX;\bC)\otimes_{\cR}\cQ'.
\end{equation}
Then
$$
\langle \bar{\one}_h, \bar{\one}_{h'}\rangle_{\cX} =
\frac{\delta_{h^{-1},h'}}{|G|\be_1}
$$
and
$$
\bar{\one}_h \star_{\cX}\bar{\one}_{h'} =\bar{\one}_{hh'}.
$$

Given $\gamma\in G^*$, define
\begin{equation}\label{eqn:bar-phi}
\bar{\phi}_\gamma := \frac{1}{|G|}\sum_{h\in G} \chi_\gamma(h^{-1}) \bar{\one}_h.
\end{equation}
Then
$$
\langle \bar{\phi}_{\gamma},\bar{\phi}_{\gamma'}\rangle_{\cX} = \frac{\delta_{\gamma\gamma'}}{|G|^2\be_1 }
$$
and
$$
\bar{\phi}_{\gamma}\star_{\cX} \bar{\phi}_{\gamma'} =  \delta_{\gamma \gamma'} \bar{\phi}_{\gamma}.
$$

In particular, $\{\bar{\phi}_\gamma:\gamma\in G^*\}$ is a canonical basis
of the semisimple $\cQ'$-algebra $H^*_{\tT,\CR}(\cX;\bC)\otimes_{\cR}\cQ'$.

\section{Generating functions}

\subsection{The total descendant Gromov-Witten potential of $\BG$}
In \cite{JK}, Jarvis-Kimura studied all genus
descendant orbifold Gromov-Witten invariants of $\BG$
for any finite group $G$. In this subsection, we recall
their results, under the assumption that $G$ is abelian.

Let $\Mbar_{g,n}(\BG)$ be moduli space of genus $g$, $n$-pointed
twisted stable maps to $\BG$, and let
$$
\ev_j:\Mbar_{g,n}(\BG)\to \IBG = \bigcup_{h\in G}(\BG)_h
$$
be the evaluation at the $j$-th marked point. Given
$h_1,\ldots, h_n\in G$, define
$$
\Mbar_{g,(h_1,\ldots, h_n)}(\BG):= \bigcap_{j=1}^n \ev_j^{-1}\left((\BG)_{h_j}\right).
$$
Then $\Mbar_{g,(h_1,\ldots, h_n)}(\BG) $ is empty unless $h_1\cdots h_n=1$.

Given $\beta_1,\ldots,\beta_n\in H^*_\CR(\BG;\bC)$
and $a_1,\ldots, a_n\in \bZ_{\geq 0}$, define the correlator
$$
\langle \tau_{a_1}(\beta_1) \cdots \tau_{a_n}(\beta_n)\rangle_{g,n}:=
\int_{\Mbar_{g,n}(\BG)}\prod_{j=1}^n \bar{\psi}_j^{a_j}\ev_j^*\beta_j.
$$
With this notation, we have
\begin{itemize}
\item For any $h_1,\ldots, h_n\in G$,
\begin{eqnarray*}
&& \langle \tau_{a_1}(\one_{h_1})\cdots \tau_{a_n}(\one_{h_n})\rangle_{g,n} \\
&=& \begin{cases}
|G|^{2g-1} \int_{\Mbar_{g,n}}\psi_1^{a_1}\cdots \psi_n^{a_n}&
\textup{if }h_1\cdots h_n=1\textup{ and } \sum_{i=1}^n a_i =3g-3+n,\\
0 &\textup{otherwise}.
\end{cases}
\end{eqnarray*}
\item For any $\gamma_1,\ldots, \gamma_n\in G^*$,
\begin{eqnarray*}
&& \langle \tau_{a_1}(\phi_{\gamma_1})\cdots \tau_{a_n}(\phi_{\gamma_n}) \rangle_{g,n} \\
&=& \begin{cases}
|G|^{2g-2}\int_{\Mbar_{g,n}}\psi_1^{a_1}\cdots \psi_n^{a_n}&
\textup{if }\gamma_1 =\cdots =\gamma_n \textup{ and } \sum_{i=1}^n a_i =3g-3+n,\\
0 &\textup{otherwise}.
\end{cases}
\end{eqnarray*}
\end{itemize}

In particular,
$$
\langle \one_h , \one_{h'}\rangle =\langle \tau_0(\one)\tau_0(\one_{h}) \tau_0(\one_{h'})\rangle_{0,3}
=\frac{\delta_{h^{-1},h'}}{|G|}.
$$
$$
\one_h \star \one_{h'} =\sum_{k\in G} \langle \one_h,\one_{h'},|G|\one_{k^{-1}}\rangle_{0,3}\one_k  =\one_{hh'}
$$
$$
\langle \phi_\gamma , \phi_{\gamma'}\rangle =\langle \tau_0(\one)\tau_0(\phi_\gamma) \tau_0(\phi_{\gamma'})\rangle_{0,3}
=\sum_{\delta\in G^*}\langle \tau_0(\phi_\delta)\tau_0(\phi_\gamma) \tau_0(\phi_{\gamma'})\rangle_{0,3}
=\frac{\delta_{\gamma,\gamma'}}{|G|^2}
$$
$$
\phi_\gamma \star \phi_{\gamma'} =\sum_{\delta\in G^*} \langle \phi_\gamma,\phi_{\gamma'},
|G|^2\phi_\delta\rangle_{0,3} \phi_\delta =\delta_{\gamma,\gamma'} \phi_\gamma.
$$

Introducing formal variables
$$
\bu=\sum_{a\geq 0}u_a z^a,
$$
where
$$
u_a =\sum_{\alpha\in G^*} u_a^\alpha \phi_\alpha \in H^*_\CR(\BG;\bC).
$$
Define
$$
\langle \bu,\ldots, \bu \rangle_{g,n}
:= \sum_{a_1,\ldots, a_n\in \bZ_{\geq 0}}  \langle\tau_{a_1}(u_{a_1})\cdots, \tau_{a_n}(u_{a_n}) \rangle_{g,n}.
$$
Then
\begin{eqnarray*}
\langle \bu,\ldots, \bu \rangle_{g,n}
&=& \sum_{a_1,\ldots, a_n\in \bZ_{\geq 0}} \sum_{\alpha_1,\ldots,\alpha_n \in G^*}
(\prod_{j=1}^n u_{a_j}^{\alpha_j}) \langle\tau_{a_1}(\phi_{\alpha_1})\cdots, \tau_{a_n}(\phi_{\alpha_n}) \rangle_{g,n} \\
&=& \sum_{\alpha\in G^*} \sum_{a_1,\ldots, a_n\in \bZ_{\geq 0}} (\prod_{j=1}^n u_{a_j}^{\alpha}) |G|^{2g-2}
\int_{\Mbar_{g,n}}\psi_1^{a_1}\cdots \psi_n^{a_n}.
\end{eqnarray*}

Define the generating functions
$$
\cF_g^{\BG}(\bu):= \sum_{n\geq 0} \frac{1}{n!} \langle \bu,\ldots, \bu\rangle_{g,n}
$$
$$
\cD^{\BG}(\bu,\hbar): =\exp\biggl(\sum_{g\geq 0}\hbar^{g-1}\cF_g^{\BG}(\bu)\biggr).
$$
The generating function $\cD^{\BG}(\bu,\hbar)$ is called the total descendant Gromov-Witten potential of $\BG$.

\subsection{The total descendant equivariant Gromov-Witten potential of $\cX=[\bC^r/G]$}
The quotient stack $\cX =[\bC^r/G]$ is the total space
of a vector bundle $F$ over $\BG$ of the form
$$
F =\bigoplus_{i=1}^r L_{\chi_i},\quad \chi_1,\ldots, \chi_r \in G^*.
$$
The vector bundle $F\longrightarrow \BG$ is equipped with a $\tT$-equivariant structure.
Let
$$
F_{g,n} := R^\bullet\pi_*f^*F=R^0\pi_*f^*F-R^1\pi_*f^*F \in K_{\tT}(\Mbar_{g,n}(\BG))
$$
be the $K$-theoretic pushforward, where $K_{\tT}$ denotes the $\tT$-equivariant $K$-theory,
$$
\pi:\cC_{g,n}\to \Mbar_{g,n}(\BG)
$$
is the projection from the universal curve and
$$
f:\cC_{g,n}\to \BG
$$
is the universal map. Similarly, define
$$
F^i_{g,n} := R^\bullet\pi_*f^*L_{\chi_i} =R^0\pi_*f^*L_{\chi_i} -R^1\pi_*f^*L_{\chi_i}
\in K_{\tT}(\Mbar_{g,n}(\BG)).
$$
Then
$$
F_{g,n}= F^1_{g,n}\oplus \cdots \oplus F^r_{g,n}.
$$
By orbifold Riemann-Roch theorem, for $i\in \{1,\ldots, r\}$,
\begin{equation}\label{eqn:rankF}
\rank_{\bC} F^i_{g,n}\Big|_{\Mbar_{g,(h_1,\ldots, h_n)}(\BG)}  = 1-g-\sum_{j=1}^n c_i(h_j).
\end{equation}

The $\tT$-equivariant Euler class of $F_{g,n}$ is
\begin{equation}\label{eqn:eT}
e_{\tT}(F_{g,n}) =\prod_{i=1}^r \left(\sw_i^{\rank F^i_{g,n}}c_{\frac{1}{\sw_i}}(F^i_{g,n})\right),
\end{equation}
where
$$
c_t(F) = 1 + c_1(F) t + c_2(F) t^2 + \cdots
$$

Given $\beta_1,\ldots,\beta_n\in H^*_{\tT,\CR}(\cX;\bC)\otimes_{\cR}{\cQ'}$
and $a_1,\ldots, a_n\in \bZ_{\geq 0}$, define $\tT$-equivariant genus $g$ descendant
orbifold Gromov-Witten invariants of $\cX$ by
\begin{equation}\label{eqn:X}
\langle \tau_{a_1}(\beta_1) \cdots \tau_{a_n}(\beta_n)\rangle^{\cX}_{g,n}:=
\int_{\Mbar_{g,n}(\BG)}\frac{\prod_{j=1}^n \bar{\psi}_j^{a_j}\ev_j^*\beta_j}{e_{\tT}(F_{g,n})}.
\end{equation}
Here we use the $\cQ'$-vector space isomorphism 
$H^*_{\tT,\CR}(\cX;\bC)\otimes_{\cR}\cQ' \cong H^*_\CR(\BG;\cQ')$ to view 
$\beta_1,\ldots,\beta_n$ as elements in $H^*_\CR(\BG;\cQ')$.
Introducing formal variables
$$
\hat{\bu}=\sum_{a\geq 0}\hat u_a z^a,
$$
where $\hat u_a \in H^*_\CR(\cX;\cQ)$, and define
$$
\langle \hat\bu,\ldots, \hat\bu \rangle^{\cX}_{g,n}
:= \sum_{a_1,\ldots, a_n\in \bZ_{\geq 0}}  \langle\tau_{a_1}(\hat u_{a_1})\cdots, \tau_{a_n}(\hat u_{a_n}) \rangle^{\cX}_{g,n}.
$$
Define the generating functions
$$
\cF_g^{\cX}(\hat\bu):= \sum_{n\geq 0} \frac{1}{n!} \langle \hat\bu,\ldots, \hat\bu\rangle^{\cX}_{g,n}
$$
$$
\cD^{\cX}(\hat\bu,\hbar): =\exp\biggl(\sum_{g\geq 0}\hbar^{g-1}\cF_g^{\cX}(\hat\bu)\biggr).
$$
The generating function $\cD^{\cX}(\hat\bu,\hbar)$ is called the total descendant
$\tT$-equivariant Gromov-Witten potential of $\cX$.

\subsection{The twisted total descendant Gromov-Witten potential}
In this subsection, we define a twisted total descendant
Gromov-Witten potential $\cD^{\tw}(\bu, \hbar)$, which is related to $D^{\cX}(\hat\bu,\hbar)$
by a simple change of variables. In Section \ref{sec:graph}, we will write
$\cD^{\tw}(\bu,\hbar)$ as a graph sum.

Given $\beta_1,\ldots, \beta_n \in H^*_\CR(\BG;\cQ')$, define
and $a_1,\ldots, a_n \in \bZ_{\geq 0}$, define twisted correlators:
\begin{equation}\label{eqn:twisted}
\langle \tau_{a_1}(\beta_1) \cdots \tau_{a_n}(\beta_n)\rangle^\tw_{g,n}:=
\int_{\Mbar_{g,n}(\BG)}\frac{\prod_{j=1}^n \bar{\psi}_j^{a_j}\ev_i^*\beta_j}{
\prod_{i=1}^r c_{\frac{1}{\sw_i}}(F^i_{g,n})}.
\end{equation}
By \eqref{eqn:rankF}, \eqref{eqn:eT}, \eqref{eqn:X}, and \eqref{eqn:twisted},
for any $h_1,\ldots, h_n\in G$,
$$
\langle \tau_{a_1}(\one_{h_1}) \cdots \tau_{a_n}(\one_{h_n})\rangle^\cX_{g,n}
=\prod_{i=1}^r \sw_i^{g-1+\sum_{j=1}^n c_i(h_j)} \langle \tau_{a_1}(\one_{h_1}) \cdots \tau_{a_n}(\one_{h_n})\rangle^\tw_{g,n}
$$
or equivalently,
\begin{equation}
\langle \tau_{a_1}(\bar{\one}_{h_1}) \cdots \tau_{a_n}(\bar{\one}_{h_n})\rangle^\cX_{g,n}
=\be_1^{g-1} \langle \tau_{a_1}(\one_1) \cdots \tau_{a_n}(\one_n)\rangle^\tw_{g,n},
\end{equation}
where $\bar{\one}_h$ is defined by \eqref{eqn:bar-one}.

Let $\phi_\gamma$ be defined as in \eqref{eqn:bar-phi}. Then
for any $\gamma_1,\ldots, \gamma_n\in G^*$,
\begin{equation}
\langle \tau_{a_1}(\bar{\phi}_{\gamma_1}) \cdots \tau_{a_n}(\bar{\phi}_{\gamma_n})\rangle^\cX_{g,n}
=\be_1^{g-1} \langle \tau_{a_1}(\phi_{\gamma_1}) \cdots \tau_{a_n}(\phi_{\gamma_n})\rangle^\tw_{g,n}.
\end{equation}

Introduce formal variables
$$
\bu=\sum_{a\geq 0}u_a z^a,
$$
where $u_a \in H^*_\CR(\cX;\cQ)$, and define
$$
\langle \bu,\ldots \bu\rangle^\tw_{g,n}
=\sum_{a_1,\ldots, a_n\in \bZ_{\geq 0}}
\langle \tau_a(u_{a_1}) \cdots \tau_{a_n}(u_{a_n})\rangle^{\tw}_{g,n}.
$$
Define the generating functions
$$
\cF_g^\tw(\bu):= \sum_{n\geq 0} \frac{1}{n!} \langle \bu,\ldots, \bu\rangle^\tw_{g,n}
$$
$$
\cD^\tw(\bu,\hbar): =\exp\biggl(\sum_{g\geq 0}\hbar^{g-1}\cF_g^\tw(\bu)\biggr).
$$
Then
$$
\cD^\cX(\hat\bu= \sum_{a\geq 0} \hat{u}_a z^a ,\hbar)
= \cD^\tw(\bu= \sum_{a\geq 0} u_a z^a,  \be_1 \hbar),
$$
where
$$
\hat{u}_a = \sum_{\alpha\in G^*} u_a^\alpha\bar{\phi}_\alpha,
\quad u_a = \sum_{\alpha\in G^*} u_a^\alpha \phi_\alpha.
$$

\section{Graph Sum} \label{sec:graph}

\subsection{Orbifold quantum Riemann-Roch}
For $i\in\{1,\ldots,r\}$, define  $\{ s^i_k: k= 1,\ldots\ \}$ by
$$
\exp\Bigl(\sum_{k\geq  1} s^i_k \frac{x^k}{k!} \Bigr) =\frac{1}{1+\frac{x}{\sw_i}}.
$$
Then
$$
s^i_k =  (k-1)!(-\sw_i)^{-k}
$$
$$
\frac{1}{\prod_{i=1}^r c_{\frac{1}{\su_i}}(F_{g,n}^i)}=\exp\Bigl(\sum_{i=1}^r \sum_{k\geq  0}  s^i_k \frac{ch_k(R^\bullet\pi_*f^*L_{\chi_i})}{k!}\Bigr).
$$

For $i\in \{1,\ldots, r\}$, and each integer $m\geq 0$, define
$$
A^i_m: H^*(\IBG;\bC)\to H^*(\IBG;\bC),\quad
\one_h \mapsto  B_m(c_i(h)) \one_h,
$$
where $B_m(x)$ is the Bernoulli polynomial. Then
$$
A^i_m(\phi_\alpha)= \sum_{\beta} (E^i_m)_\alpha^\beta \phi_\beta,
$$
where
$$
(E^i_m)_\alpha^\beta =\frac{1}{|G|}\sum_{h\in G}\chi_\alpha(h^{-1})B_m(c_i(h)) \chi_\beta(h).
$$

For each $i\in \{1,\ldots,r\}$,  $A^i_0, A^i_2, A^i_4,\ldots$ are self-adjoint, and
$A^i_3, A^i_5, A^i_7,\ldots$ are anti-self-adjoint.
By Tseng's orbifold quantum Riemann-Roch theorem \cite{Ts10},
\begin{equation}\label{eqn:QRR}
\cD^\tw = \exp\left(\sum_{m\geq 1} \sum_{i=1}^r \frac{s^i_m}{(m+1)!}(A^i_{m+1}z^m)^{\wedge}\right) \cD^{\BG},
\end{equation}
where
\begin{eqnarray*}
&& \bigl(A^i_{m+1}z^m\bigr)^{\wedge} \\
&=& \sum_{\alpha,\beta\in G^*} (E^i_{m+1})^\alpha_\beta  \left( -\sum_{\ell\geq 0}
q^\beta_l \frac{\partial}{\partial q^\alpha_{\ell +m}}
+\frac{\hbar|G|^2}{2}\sum_{\ell=0}^{m-1} \sum_{\alpha,\beta\in G^*} (-1)^{\ell+1+m}
\frac{\partial}{\partial q^\alpha_\ell}\frac{\partial}{\partial q^\beta_{m-1-\ell} }\right).
\end{eqnarray*}
We have the dilaton shift
$$
\sum_{\beta \in G^*} q_1^\beta\phi_\beta =\sum_{\beta \in G^*}u_1^\beta\phi_\beta -\one
=\sum_{\beta \in G^*}(u_1^\beta-1)\phi_\beta.
$$
So
$$
q_1^\beta = u_1^\beta-1,\quad \forall \beta \in G^*.
$$
Therefore, we have
\begin{eqnarray*}
\bigl(A^i_{m+1}z^m\bigr)^{\wedge} &=&
\sum_{\alpha,\beta\in G^*} (E^i_{m+1})^\alpha_\beta  \left(\frac{\partial}{\partial u^\alpha_m}-  \sum_{\ell\geq 0}
u^\beta_\ell \frac{\partial}{\partial u^\alpha_{\ell +m}} \right.\\
&& \quad\quad + \left.\frac{\hbar|G|^2}{2}\sum_{\ell=0}^{m-1} \sum_{\alpha,\beta\in G^*} (-1)^{\ell+1+m}
\frac{\partial}{\partial u^\alpha_\ell}\frac{\partial}{\partial u^\beta_{m-1-\ell} }\right).
\end{eqnarray*}

\subsection{Sum over graphs} \label{sec:graph-sum}
In this section, we write the right hand side of \eqref{eqn:QRR} as a sum
over graphs, following \cite{DSS, DOSS}.

Given a connected graph $\Ga$, we introduce the following notation.
\begin{enumerate}
\item $V(\Ga)$ is the set of vertices in $\Ga$.
\item $E(\Ga)$ is the set of edges in $\Ga$.
\item $H(\Ga)$ is the set of half edges in $\Gamma$.
\item $L^o(\Ga)$ is the set of ordinary leaves in $\Ga$.
\item $L^1(\Ga)$ is the set of dilaton leaves in $\Ga$.
\end{enumerate}

With the above notation, we introduce the following labels:
\begin{enumerate}
\item (genus) $g: V(\Ga)\to \bZ_{\geq 0}$.
\item (marking) $\alpha: V(\Ga) \to G^*$. This induces
$\alpha:L(\Ga)=L^o(\Ga)\cup L^1(\Ga)\to G^*$, as follows:
if $l\in L(\Ga)$ is a leaf attached to a vertex $v\in V(\Ga)$,
define $\alpha(l)=\alpha(v)$.
\item (height) $k: H(\Ga)\to \bZ_{\geq 0}$.
\end{enumerate}

Given an edge $e$, let $h_1(e),h_2(e)$ be the two half edges associated to $e$. The order of the two half edges does not affect the graph sum formula in this paper. Given a vertex $v\in V(\Ga)$, let $H(v)$ denote the set of half edges
emanating from $v$. The valency of the vertex $v$ is equal to
the size of the set $H(v)$: $\val(v)=|H(v)|$.
A labeled graph $\vGa=(\Ga,g,\alpha,k)$ is {\em stable} if
$$
2g(v)-2 + \val(v) >0
$$
for all $v\in V(\Ga)$.

Let $\bGa(\BG)$ denote the set of all stable labeled graphs
$\vGa=(\Gamma,g,\alpha,k)$. The genus of a stable labeled graph
$\vGa$ is defined to be
$$
g(\vGa):= \sum_{v\in V(\Ga)}g(v)  + |E(\Ga)|- |V(\Ga)|  +1
=\sum_{v\in V(\Ga)} (g(v)-1) + (\sum_{e\in E(\Gamma)} 1) +1.
$$
Define
$$
\bGa_{g,n}(\BG)=\{ \vGa=(\Gamma,g,\alpha,k)\in \bGa(\BG): g(\vGa)=g, |L^o(\Ga)|=n\}.
$$

Define the operator $R(z)$ to be
$$
R(z)=   \exp\Bigl(\sum_{m\geq 1}\frac{(-1)^{m}}{m(m+1)} \sum_{i=1}^r A^i_{m+1} (\frac{z}{\sw_i})^m \Bigr).
$$
Then \eqref{eqn:QRR} can be written as
$$
\cD^\tw =R(z)^\wedge  \, \cD^{\BG}.
$$
Under the canonical basis $\{\phi_\gamma:\gamma\in G^*\}$, we have
$$
R(z)^\alpha_\beta =\frac{1}{|G|}\sum_{h\in G} \chi_\alpha(h)
\chi_\beta(h^{-1})\exp\Bigl(\sum_{m=1}^\infty \frac{(-1)^m}{m(m+1)}\sum_{i=1}^r B_{m+1}(c_i(h))(\frac{z}{\sw_i})^m \Bigr).
$$
Given $\beta\in G^*$, define
$$
\bu^\beta(z) = \sum_{a\geq 0} u^\beta_a z^a.
$$

We assign weights to leaves, edges, and vertices of a labeled graph $\vGa\in \bGa(\BG)$ as follows.
\begin{enumerate}
\item {\em Ordinary leaves.} To each ordinary leaf $l \in L^o(\Ga)$ with  $\alpha(l)= \alpha\in G^*$
and  $k(l)= k\in \bZ_{\geq 0}$, we assign:
$$
(\tcL^{\bu})^\alpha_k(l) = [z^k] (\sum_{\beta\in G^*}R(-z)^\beta_\alpha \bu^\beta(z)).
$$
\item {\em Dilaton leaves.} To each dilaton leaf $l \in L^1(\Ga)$ with $\alpha(l)=\alpha \in G^*$
and $2\leq k(l)=k \in \bZ_{\geq 0}$, we assign
$$
(\tcL^1)^\alpha_k(l) = [z^{k-1}](-\sum_{\beta\in G^*} R(-z)^\beta_\alpha).
$$
Note that
$$
\frac{1}{|G|}\sum_{\beta\in G^*}\chi_\beta(h) = \frac{1}{|G|}\sum_{\beta\in G^*} \chi_\beta(h)\chi_\beta(1)=\delta_{h,1}.
$$
So
$$
\sum_{\beta\in G^*}R^\beta_\alpha(-z) = 
\exp\Bigl(\sum_{m=1}^\infty \frac{1}{m(m+1)}\sum_{i=1}^r B_{m+1}(\frac{z}{\sw_i })^m \Bigr).
$$

\item {\em Edges.} To an edge connected a vertex marked by $\alpha\in G^*$ to a vertex
marked by $\beta\in G^*$ and with heights $k$ and $l$ at the corresponding half-edges, we assign
$$
\widetilde{\cE}^{\alpha,\beta}_{k,l}(e) = [z^k w^l]
\Bigl(\frac{|G|^2}{z+w} (\delta_{\alpha,\beta}-\sum_{\gamma\in G^*} R(-z)^\gamma_\alpha R(-w)^\gamma_\beta)\Bigr).
$$
\item {\em Vertices.} To a vertex $v$ with genus $g(v)=g\in \bZ_{\geq 0}$ and with
marking $\alpha(v)=\gamma\in G^*$, with $n$ ordinary
leaves and half-edges attached to it with heights $k_1, ..., k_n \in \bZ_{\geq 0}$ and $m$ more
dilaton leaves with heights $k_{n+1}, \ldots, k_{n+m}\in \bZ_{\geq 0}$, we assign
$$
|G|^{2g-2} \int_{\Mbar_{g,n+m}}\psi_1^{k_1} \cdots \psi_{n+m}^{k_{n+m}}.
$$
\end{enumerate}

We define the weight of a labeled graph $\vGa\in \bGa(\BG)$ to be
\begin{eqnarray*}
\widetilde{w}(\vGa) &=& \prod_{v\in V(\Ga)} |G|^{2g-2} \langle \prod_{h\in H(v)} \tau_{k(h)}\rangle_{g(v)}
\prod_{e\in E(\Ga)} \widetilde{\cE}^{\alpha(v_1(e)),\alpha(v_2(e))}_{k(h_1(e)),k(h_2(e))}(e)\\
&& \cdot \prod_{l\in L^o(\Ga)}(\tcL^{\bu})^{\alpha(l)}_{k(l)}(l)
\prod_{l\in L^1(\Ga)}(\tcL^1)^{\alpha(l)}_{k(l)}(l).
\end{eqnarray*}
Then
$$
\log(\cD^\tw(\bu))=\log(\hat{R} \cD^{\BG}(\bu)) =
\sum_{\vGa\in \bGa(\BG)}  \frac{ \hbar^{g(\vGa)-1} \widetilde{w}(\vGa)}{|\Aut(\vGa)|}
= \sum_{g\geq 0}\hbar^{g-1} \sum_{n\geq 0}\sum_{\vGa\in \bGa_{g,n}(\cB G)}\frac{\widetilde{w}(\vGa)}{|\Aut(\vGa)|}.
$$
$$
\log(\cD^\cX(\hat{\bu}))=
\sum_{\vGa\in \bGa(\BG)}  \frac{ (\hbar\be_1)^{g(\vGa)-1} \widetilde{w}(\vGa)}{|\Aut(\vGa)|}.
$$
To obtain the graph sum of $\log(\cD^\cX(\hat\bu))$, we assign a new weight $w(\vGa)$
to each labeled graph $\vGa\in \bGa(\BG)$.  Define
\begin{eqnarray*}
&& (\cL^\bu)^\alpha_k(l) = \frac{1}{|G|\sqrt{\be_1}} (\tcL^\bu)^\alpha_k(l), \quad
(\cL^1)^\alpha_k(l)  = \frac{1}{|G|\sqrt{\be_1}} (\tcL^1)^\alpha_k(l), \\
&& \cE^{\alpha,\beta}_{k,l}(e) = [z^k w^l]
\Bigl(\frac{1}{z+w} (\delta_{\alpha,\beta}-\sum_{\gamma\in G^*} R(-z)^\gamma_\alpha R(-w)^\gamma_\beta)\Bigr),
\end{eqnarray*}
and define
\begin{eqnarray*}
w(\vGa) &=& \prod_{v\in V(\Ga)} (|G|\sqrt{\be_1})^{2g(v)-2+\val(v)} \langle \prod_{h\in H(v)} \tau_{k(h)}\rangle_{g(v)}
\prod_{e\in E(\Ga)}\cE^{\alpha(v_1(e)),\alpha(v_2(e))}_{k(h_1(e)),k(h_2(e))}(e)\\
&& \cdot \prod_{l\in L^o(\Ga)}(\cL^{\bu})^{\alpha(l)}_{k(l)}(l)
\prod_{l\in L^1(\Ga)}(\cL^1)^{\alpha(l)}_{k(l)}(l).
\end{eqnarray*}
With the above notation, we have the following graph sum for $\log(\cD^{\cX}(\hat\bu))$.
\begin{theorem}\label{main}
\begin{equation}
\log(\cD^\cX(\hat\bu))
= \sum_{g\geq 0}\hbar^{g-1} \sum_{n\geq 0}\sum_{\vGa\in \bGa_{g,n}(\BG)}\frac{w(\vGa)}{|\Aut(\vGa)|}.
\end{equation}
\begin{equation} \label{eqn:n-point}
\frac{\langle \hat\bu,\ldots,\hat\bu\rangle_{g,n}^{\cX}}{n!} =\sum_{\vGa\in\bGa_{g,n}(\BG)}\frac{w(\vGa)}{|\Aut(\vGa)|}.
\end{equation}
\end{theorem}

Finally, we state a straightforward generalization of \eqref{eqn:n-point}. Let
$\bGa_{g,n,n'}(\BG)$ be the set of genus $g$ stable graphs
with $n$ {\em ordered} ordinary leaves and $n'$ {\em unordered} ordinary leaves.
In particular, $\bGa_{g,0,n}(\BG)=\bGa_{g,n}(\BG)$.
Given $\vGa\in \bGa_{g,n,n'}(\BG)$, let $L^O(\Gamma)=\{l_1\ldots, l_n\}$
be the set of ordered ordinary leaves in $\Gamma$, and
let $L^o(\Gamma)$ be the set of unordered ordinary leaves in $\Gamma$;
we assign the following weight to this labeled graph $\vGa\in \Gamma_{g,n,n'}(\BG)$:
\begin{eqnarray*}
w(\vGa) &=& \prod_{v\in V(\Ga)} (|G|\sqrt{\be_1})^{2g(v)-2+\val(v)} \langle \prod_{h\in H(v)} \tau_{k(h)}\rangle_{g(v)}
\prod_{e\in E(\Ga)}\cE^{\alpha(v_1(e)),\alpha(v_2(e))}_{k(h_1(e)),k(h_2(e))}(e)\\
&& \cdot \prod_{j=1}^n(\cL^{\bu_j})^{\alpha(l_j)}_{k(l_j)}(l_j)
\prod_{l\in L^o(\Ga)}(\cL^{\bu})^{\alpha(l)}_{k(l)}(l)
\prod_{l\in L^1(\Ga)}(\cL^1)^{\alpha(l)}_{k(l)}(l),
\end{eqnarray*}
where
$$
(\cL^{\bu_j})^\alpha_k(l_j) = [z^k]
\frac{1}{|G|\sqrt{\be_1}}\Big(\sum_{\beta\in G^*} R(-z)^\beta_\alpha \bu_j^\beta(z)\Big),
\quad
\bu^\beta_j (z) = \sum_{a\geq 0} (u_j)^\beta_a z^a.
$$
Define
$$
\hbu_j =\sum_{\beta\in G^*}\bu^\beta_j(z) \bar{\phi}_\beta.
$$
We have the following generalization of \eqref{eqn:n-point}.
\begin{theorem} \label{ordered}
$$
\frac{1}{(n')!}\langle \hbu_1,\ldots,\hbu_n,\bu,\ldots,\bu \rangle_{g,n+n'}^{\cX} =\sum_{\vGa\in\bGa_{g,n,n'}(\BG)}\frac{w(\vGa)}{|\Aut(\vGa)|}.
$$
\end{theorem}
Theorem \ref{ordered} will play a crucial role in the proof of the remodeling conjecture
for affine toric Calabi-Yau 3-orbifolds in \cite{FLZ}.

\end{document}